\documentclass[12pt]{iopart}
\usepackage{wrapfig}

\begin{document}

\title[Linearization of fourth-order ODEs  ]{Linearization of
fourth-order ordinary differential equations by point
transformations }

\author{Nail H. Ibragimov$^1$, Sergey V. Meleshko$^2$ and \\
Supaporn Suksern$^2$}

\address{$^1$ Mathematics and Science, Research Centre ALGA: Advances in Lie Group
Analysis, Blekinge Institute of Technology, SE-371 79 Karlskrona,
Sweden}
\address{$^2$ School of Matematics, Suranaree University of
Technology, Nakhon Ratchasima, 30000, Thailand}

\begin{abstract}
 We present here the solution of the problem on linearization of fourth-order
 equations by means of point transformations.
 We show that all fourth-order equations that are linearizable by point
 transformations are contained  in the class of equations
 which is linear in the third-order derivative.  We provide the linearization
 test and describe the procedure for obtaining the
 linearizing transformations as well as the linearized equation.
\end{abstract}


\maketitle
\section{Introduction}
The problem on linearization of second-order ordinary differential
equations by means of point transformations was solved by Sophus Lie
\cite{Lie} in 1883. More specifically, he showed that the
linearizable equations are at most cubic in the first-order
derivative and gave the linearization test in terms of the
coefficients of these equations.

In 1997, G. Grebot \cite{Gre} studied the linearization of
third-order equations by means of a restricted class of point
transformations, namely $ t = \varphi \left( x \right), u = \psi
\left( {x,y} \right). $ However, the problem was not completely
solved.

In 2004, N.H. Ibragimov and S.V. Meleshko \cite{Ibra-Mel} solved the
problem of linearization of third-order equations by means of point
transformations.  They showed that all third-order equations that
are linearizable by point transformations are contained either in
the class of equations which is linear in the second-order
derivative, or in the class of equations which is quadratic in the
second-order derivative. They provided the linearization test for
each of these classes and describe the procedure for obtaining the
linearizing transformations as well as the linearized equation.

The present paper is devoted to obtain criteria for a fourth-order
equation to be linearizable by change of the dependent and
independent variables. In our calculations we used computer algebra
packages. The final results were checked by comparing with
theoretical results on invariants as well as by applying to numerous
known and new examples of linearization. The paper is organized as
follows.

\section{Point transformations of fourth-order equations}
We consider the fourth-order ordinary differential equation
 \begin{equation}
 \label{eq:01}
    y^{\left( 4 \right)}  = f\left( {x,y,y',y'',y'''} \right).
 \end{equation}
We apply a point transformation
\begin{equation}
 \label{eq:02}
t = \varphi \left( {x,y} \right),\quad u = \psi \left( {x,y} \right)
 \end{equation}
 to equation (\ref{eq:01}).

We begin with investigating the necessary conditions for
linearization. The general form of  (\ref{eq:01}) that can be
obtained from linear equations by any point transformations
(\ref{eq:02}) is found on this step. In consequence, we identify two
candidates for linearization.

A linear fourth-order ordinary differential equation we use in the
Laguerre form.  In 1879, E. Laguerre showed that in linear ordinary
differential equation of order $n \geq 3$ the two terms of orders
next below the highest can be simultaneously removed by equivalence
transformation (see \cite{Ibra}, Section 10.2.1 and the references
therein). Therefore, we write the general linear fourth-order
equation in Laguerre's form
 \begin{equation}
 \label{eq:03}
u^{\left( 4 \right)}  + \alpha \left( t \right)u' + \beta \left( t
\right)u = 0,
 \end{equation}
  where $t$ and $u$ are  the independent and dependent variables,
  respectively.

\subsection{The candidates for linearization}
 \label{candidates}

 Considering $t$ and $u$ as the new independent and dependent
variables, respectively, one obtains the following transformation of
the first-order derivative
 \begin{equation}
 \label{eq:04}
 u' = \frac{D_x(\psi)}{D_x(\varphi)} =
 \frac{\psi_x + y'\psi_y}{\varphi_x + y'\varphi_y}\,,
 \end{equation}
 where
 $\varphi_x = \frac{\partial \varphi}{\partial x}, ~\varphi_y = \frac{\partial \varphi}{\partial
 y},$ etc., and
\[
 D_x = \frac{\partial}{\partial x} + y'\frac{\partial}{\partial y}
  + y'' \frac{\partial}{\partial y'}  + y''' \frac{\partial}{\partial y''}
  + y^{(4)}\frac{\partial}{\partial y'''}+\cdots
 \]
 is the total derivative.
 Likewise, one obtains the transformation of derivatives of the second and higher order.
 Namely, denoting by $P (x, y, y')$ the right-hand side of
 (\ref{eq:04}),
 \[
 P (x, y, y') = \frac{\psi_x + y'\psi_y}{\varphi_x + y'\varphi_y}\,
 \]
 one has
\begin{equation}
 \label{eq:05}
 u'' = \frac{{D_x \left( P \right)}} {{D_x
\left( \varphi  \right)}} = \frac{P_{x} + y'P_{y} +
y''P_{y'}}{\varphi_x + y'\varphi_y} = \frac{\Delta } {{\left(
{\varphi _x + y'\varphi _y } \right)^3 }}y'' +  \cdots.
 \end{equation}
 Denoting by $Q (x, y, y', y'')$ the right-hand side of
 (\ref{eq:05}),
 \[
 Q = \frac{\Delta}{(\varphi_x + y'
 \varphi_y)^3}\,y'' + \cdots\,
 \]
 one has
 \begin{equation}
 \label{eq:06}
 \begin{array}{ll}
  u''' & = \frac{D_x(Q)}{D_x(\varphi)} = \frac{Q_{x} + y'Q_{y} +
y''Q_{y'} + y'''Q_{y''}}{\varphi_x + y'\varphi_y}\\& =
 \frac{\Delta}{(\varphi_x + y'
 \varphi_y)^5}\,\Big[ (\varphi_x + y' \varphi_y)\,y''' - 3 \varphi_y\,(y'')^2\Big] +
 \cdots.
 \end{array}
 \end{equation}
 Denoting by $R (x, y, y', y'', y''')$ the right-hand side of
(\ref{eq:06}),
  \[
R = \frac{\Delta}{(\varphi_x + y'
 \varphi_y)^5}\,\Big[ (\varphi_x + y' \varphi_y)\,y''' - 3 \varphi_y\,(y'')^2\Big] +
 \cdots
\]
hence,
\begin{equation}
\label{eq:07}
\begin{array}{ll}
 u^{(4)} & = \frac{D_x(R)}{D_x(\varphi)} = \frac{R_{x}
+ y'R_{y} + y''R_{y'} + y'''R_{y''} + y^{(4)}R_{y'''}}{\varphi_x +
y'\varphi_y} \\&=  \frac{\Delta}{(\varphi_x + y'
 \varphi_y)^7}\,\Big[   ( \varphi_{x} + y'\varphi_{y} )^2 y^{(4)}
\Big] +
 \cdots.
\end{array}
\end{equation}
Thus, (\ref{eq:03}) becomes
\begin{equation}
\eqalign{\frac{1}{(\varphi_x + y'
 \varphi_y)^7}\Big[(\varphi_{x} + \varphi_{y}y')^2\Delta y^{(4)} +[{-
10\Delta(\varphi_{x} + \varphi_{y}y')\varphi_{y}}y'' \cr {-
2(2(5\varphi_{xy}\Delta - \varphi_{y}\Delta_{x})\varphi_{y} +
(5\varphi_{yy}\Delta - 4\varphi_{y} \Delta_{y})\varphi_{x})}y'^2 \cr
{- 2(5(2\varphi_{xy}\varphi_{x} + \varphi_{xx}\varphi_{y})\Delta -
2(\varphi_{x}\Delta_{y} + 2\varphi_{y} \Delta_{x})\varphi_{x})}y'
\cr
+ \cdots]y''' + \cdots \Big]=0.}\label{eq:00}
\end{equation}
Here
\[
\Delta = \varphi_x \psi_y - \varphi_y \psi_x \neq 0
\]
is the Jacobian of the change of variables (\ref{eq:02}). It is
manifest from  (\ref{eq:00}) that the transformations (\ref{eq:02})
with $\varphi_y= 0$ and $\varphi_y \not= 0,$ respectively, provide
two distinctly different candidates for linearization.

If $\varphi_y = 0$ we work out the missing terms in (\ref{eq:00}),
substitute the resulting expression in (\ref{eq:03})  and obtain the
following equation
 \begin{equation}
 \label{eq:09}
 \begin{array}{ll}
  y^{(4)}  & + (A_{1}y' + A_{0})y''' + B_{0}y''^2 + (C_{2}y'^2 + C_{1}y' + C_{0})y''
 \\ &+ D_{4}y'^4 + D_{3}y'^3 + D_{2}y'^2 + D_{1}y' + D_{0} = 0,
\end{array}
 \end{equation}
 where
 \begin{eqnarray}
 \fl A_{1} =  4(\psi_{y})^{-1}\psi_{yy},\label{eq:10}\\[1.5ex]
 \fl A_{0} =  - 2(\varphi_{x}\psi_{y})^{-1}(3\varphi_{xx}\psi_{y} -
2\varphi_{x}\psi_{xy}),\label{eq:11}\\[1.5ex]
 \fl B_{0} =  3(\psi_{y})^{-1}\psi_{yy},\label{eq:12}\\[1.5ex]
 \fl C_{2} =  6(\psi_{y})^{-1}\psi_{yyy},\label{eq:13}\\[1.5ex]
 \fl C_{1} =  - 6(\varphi_{x}\psi_{y})^{-1}(3\varphi_{xx}\psi_{yy} -
2\varphi_{x}\psi_{xyy}),\label{eq:14}\\[1.5ex]
\fl C_{0} =  -
(\varphi_{x}^2\psi_{y})^{-1}\Big[(4\varphi_{xxx}\varphi_{x} -
15\varphi_{xx}^2)\psi_{y} + 6(3\varphi_{xx}\psi_{xy} -
\varphi_{x}\psi_{xxy})\varphi_{x}\Big],\label{eq:15}\\[1.5ex]
\fl D_{4} =  (\psi_{y})^{-1}\psi_{yyyy},\label{eq:16}\\[1.5ex]
\fl D_{3} = - 2(\varphi_{x}\psi_{y})^{-1}(3\varphi_{xx}\psi_{yyy} -
2\varphi_{x}\psi_{xyyy}),\label{eq:17}\\[1.5ex]
\fl D_{2} =  -
(\varphi_{x}^2\psi_{y})^{-1}(4\varphi_{xxx}\varphi_{x}\psi_{yy} -
15\varphi_{xx}^2\psi_{yy} + 18\varphi_{xx}\varphi_{x}\psi_{xyy} -
6\varphi_{x}^2\psi_{xxyy}),\label{eq:18}\\[1.5ex]
\fl  D_{1}   =  -
(\varphi_{x}^3\psi_{y})^{-1}\Big[3(5\varphi_{xx}^2\psi_{y} -
10\varphi_{xx}\varphi_{x}\psi_{xy} +
6\varphi_{x}^2\psi_{xxy})\varphi_{xx}
 -(\varphi_{x}^3\psi_{y}\alpha + 4\psi_{xxxy})\varphi_{x}^3 \nonumber \\-
2(5\varphi_{xx}\psi_{y} -
4\varphi_{x}\psi_{xy})\varphi_{xxx}\varphi_{x} +
\varphi_{xxxx}\varphi_{x}^2\psi_{y}\Big],\label{eq:19}\\
\fl   D_{0}  =  - (\varphi_{x}^3\psi_{y})^{-1}\Big[(15\varphi_{xx}^3
- \varphi_{x}^6\alpha + \varphi_{xxxx}\varphi_{x}^2)\psi_{x}  -
(10\varphi_{xxx}\varphi_{xx}\psi_{x} -
4\varphi_{xxx}\varphi_{x}\psi_{xx} \nonumber \\+
15\varphi_{xx}^2\psi_{xx} - 6\varphi_{xx}\varphi_{x}\psi_{xxx} +
\varphi_{x}^6\beta\psi +
\varphi_{x}^2\psi_{xxxx})\varphi_{x}\Big].\label{eq:20}
\end{eqnarray}
{\bf Definition 1}.
 We call  (\ref{eq:09}) with arbitrary coefficients
 $A_0 = A_0(x, y),$\linebreak $ A_1 = A_1(x, y), \;\;
B_0 = B_0(x, y),\;\; C_0 = C_0(x, y),\;\; C_1 = C_1(x, y),\;\; C_2 =
C_2(x, y),$ and
 $D_i = D_i(x, y), (i = 0, \ldots, 4),$ the
 first candidate for linearization.

If $\varphi_y \not= 0,$  we proceed likewise and setting $r(x, y)
=\frac{\varphi_x}{\varphi_y} $, arrive at the following equation
\begin{equation}
 \label{eq:21}
\begin{array}{c}
  y^{( 4)}   +  \frac{1}{y'  +  r}
(  - 10 y''  +  F_2 y'^2  +  F_1 y'  +  F_{0 } )y'''  \hfill \\
   + \frac{1}
{{( y'  +  r)^2 }}\left[
  15y''^3   +  ( H_2 y'^2  +  H_1 y' +  H_0 )y''^2  \right.\hfill \\
   +  (J_4 y'^4  +  J_3 y'^3  +  J_2 y'^2  +  J_1 y'  + J_{0 } )y'' \hfill \\
   +  K_7 y'^7  +  K_6 y'^6  +  K_5 y'^5  +  K_4 y'^4  \hfill \\
\left.   +  K_3 y'^3  +  K_2 y'^2  +  K_1 y' +  K_0
  \right] = 0, \hfill \\
\end{array}
\end{equation}
where
\begin{eqnarray}
  \fl  F_{2} =  - 2(\varphi_{y}\Delta)^{-1}(5\varphi_{yy}\Delta
 - 2\varphi_{y}\Delta_{y}), \label{eq:23}\\[1.5ex]
  \fl  F_{1} =  4(\varphi_{y}\Delta)^{-1}\Big[(\Delta_{x} + \Delta_{y}r
-5r_{y}\Delta)\varphi_{y} - 5\varphi_{yy}r\Delta\Big], \label{eq:24}\\[1.5ex]
 \fl   F_{0} =   - 2(\varphi_{y}\Delta)^{-1}\Big[((5r_{y}\Delta -
2\Delta_{x})r + 5r_{x}\Delta)\varphi_{y} + 5\varphi_{yy}r^2\Delta\Big], \label{eq:25}\\[1.5ex]
  \fl  H_{2} =  6(\varphi_{y}\Delta)^{-1}(5\varphi_{yy}\Delta -
2\varphi_{y}\Delta_{y}), \label{eq:27}\\[1.5ex]
 \fl   H_{1} =  - 3(\varphi_{y}\Delta)^{-1}\Big[(5\Delta_{x} +
3\Delta_{y}r
- 25r_{y}\Delta)\varphi_{y} - 20\varphi_{yy}r\Delta\Big], \label{eq:28}\\[1.5ex]
 \fl   H_{0} =  3(\varphi_{y}\Delta)^{-1}\Big[(5(3r_{x} + 2r_{y}r)\Delta
- (5\Delta_{x} - \Delta_{y}r)r)\varphi_{y} +
10\varphi_{yy}r^2\Delta\Big], \label{eq:29}\\[1.5ex]
\fl    J_{4} =  -
(\varphi_{y}^2\Delta)^{-1}(10\varphi_{yyy}\varphi_{y}\Delta -
45\varphi_{yy}^2\Delta + 30\varphi_{yy}\varphi_{y}\Delta_{y} -
6\varphi_{y}^2\Delta_{yy}), \label{eq:30}\\[1.5ex]
 \fl    J_{3} =
2(\varphi_{y}^2\Delta)^{-1}\Big[3((2(\Delta_{xy} + \Delta_{yy}r -
5r_{y}\Delta_{y}) - 5r_{yy}\Delta)\varphi_{y}^2 \nonumber
\\  - 5((\Delta_{x} + 3\Delta_{y}r - 4r_{y}\Delta)\varphi_{y} -
6\varphi_{yy}r\Delta)\varphi_{yy})
 - 20\varphi_{yyy}\varphi_{y}r\Delta\Big], \label{eq:31}\\[1.5ex]
 \fl   J_{2} =  6(\varphi_{y}^2\Delta)^{-1}\Big[(\Delta_{xx} +
\Delta_{yy}r^2 + 4\Delta_{xy}r - 5(2\Delta_{x} + 3\Delta_{y}r -
5r_{y}\Delta)r_{y}\nonumber \\
 - 10r_{yy}r\Delta - 5r_{x}\Delta_{y} - 5r_{xy}\Delta)\varphi_{y}^2
- 5(((3(\Delta_{x} + \Delta_{y}r) - 10r_{y}\Delta)r \nonumber\\  -
2r_{x}\Delta)\varphi_{y} - 9\varphi_{yy}r^2\Delta)\varphi_{yy} -
10\varphi_{yyy}\varphi_{y}r^2\Delta\Big], \label{eq:32}\\[1.5ex]
\fl J_{1} = - 2(\varphi_{y}^2\Delta)^{-1}\Big[((5(3(3\Delta_{x} +
\Delta_{y}r) - 14r_{y}\Delta)r_{y} - 6(\Delta_{xy}r + \Delta_{xx})
\nonumber \\  + 20r_{yy}r\Delta)r + 5(3(\Delta_{x} + \Delta_{y}r) -
16r_{y}\Delta)r_{x} + 5r_{xx}\Delta + 20r_{xy}r\Delta)\varphi_{y}^2
\nonumber \\  + 15(((3\Delta_{x} + \Delta_{y}r - 8r_{y}\Delta)r -
4r_{x}\Delta)\varphi_{y} - 6\varphi_{yy}r^2\Delta)\varphi_{yy}r
\nonumber \\  + 20\varphi_{yyy}\varphi_{y}r^3\Delta\Big],
\label{eq:33}
\end{eqnarray}
\begin{eqnarray}\label{eq:34}
\fl J_{0} =   - (\varphi_{y}^2\Delta)^{-1}\Big[((2((5r_{yy}r\Delta -
3\Delta_{xx})r + 5r_{xx}\Delta + 5r_{xy}r\Delta) \nonumber \\  -
5(7r_{y}\Delta - 6\Delta_{x})r_{y}r)r - 5(2(7r_{y}\Delta
 - 3\Delta_{x})r + 9r_{x}\Delta)r_{x})\varphi_{y}^2
 \nonumber \\  - 5(3(2((2r_{y}\Delta - \Delta_{x})r + 2r_{x}\Delta)\varphi_{y}
 + 3\varphi_{yy}r^2\Delta)\varphi_{yy}
 \nonumber
 \\
 - 2\varphi_{yyy}\varphi_{y}r^2\Delta)r^2\Big],
\end{eqnarray}
\begin{eqnarray}\label{eq:35}
\fl  K_{7} =  -
(\varphi_{y}^2\Delta)^{-1}\Big[\varphi_{yyyy}\varphi_{y}^2\psi_{y} -
10\varphi_{yyy}\varphi_{yy}\varphi_{y}\psi_{y} +
4\varphi_{yyy}\varphi_{y}^2\psi_{yy} + 15\varphi_{yy}^3\psi_{y}
\nonumber \\  - 15\varphi_{yy}^2\varphi_{y}\psi_{yy} +
6\varphi_{yy}\varphi_{y}^2\psi_{yyy} - \varphi_{y}^7\beta\psi -
\varphi_{y}^6\psi_{y}\alpha - \varphi_{y}^3\psi_{yyyy}\Big],
\end{eqnarray}
\begin{eqnarray} \label{eq:36}
\fl  K_{6} =
(\varphi_{y}^3\Delta)^{-1}\Big[3(5((7\varphi_{y}\psi_{yy}r -
6\Delta_{y})\varphi_{y} - 7(\varphi_{y}\psi_{y}r -
\Delta)\varphi_{yy})\varphi_{yy} \nonumber \\  -
2(7\varphi_{y}\psi_{yyy}r - 5\Delta_{yy})\varphi_{y}^2)\varphi_{yy}
+ (7\varphi_{y}^5 \beta \psi r + 7\varphi_{y}^4\psi_{y}\alpha r  -
\varphi_{y}^3\alpha\Delta \nonumber \\
  + 7\varphi_{y}\psi_{yyyy}r  -
4\Delta_{yyy})\varphi_{y}^3 + 2(35\varphi_{yy}\varphi_{y}\psi_{y}r
 - 30\varphi_{yy}\Delta - 14\varphi_{y}^2\psi_{yy}r \nonumber \\
  + 10\varphi_{y}\Delta_{y})\varphi_{yyy}\varphi_{y}
 - (7\varphi_{y}\psi_{y}r -
 5\Delta)\varphi_{yyyy}\varphi_{y}^2\Big]
,
\end{eqnarray}
\begin{eqnarray}\label{eq:37}
 \fl K_{5} = -
(\varphi_{y}^3\Delta)^{-1}\Big[(2(3(\Delta_{xyy} + 3\Delta_{yyy}r -
5r_{y}\Delta_{yy} - 5r_{yy}\Delta_{y}) -
5r_{yyy}\Delta) \nonumber \\
  - 3(7\varphi_{y}^4\beta\psi r    +
7\varphi_{y}^3\psi_{y}\alpha r - 2\varphi_{y}^2\alpha\Delta +
7\psi_{yyyy}r)\varphi_{y}r)\varphi_{y}^3 \nonumber \\
  - 3(2(5(\Delta_{xy} +
5\Delta_{yy}r - 4r_{y}\Delta_{y} - 2r_{yy}\Delta)  -
21\varphi_{y}\psi_{yyy}r^2)\varphi_{y}^2 \nonumber \\
  - 15((\Delta_{x} +
11\Delta_{y}r - 3r_{y}\Delta   -
7\varphi_{y}\psi_{yy}r^2)\varphi_{y} \nonumber \\
 + 7(\varphi_{y}\psi_{y}r -
2\Delta)\varphi_{yy}r)\varphi_{yy})\varphi_{yy}
 - 2((5(\Delta_{x}
+ 11\Delta_{y}r - 3r_{y}\Delta) \nonumber \\
  -
42\varphi_{y}\psi_{yy}r^2)\varphi_{y} + 15(7\varphi_{y}\psi_{y}r
 - 12\Delta)\varphi_{yy}r)\varphi_{yyy}\varphi_{y}
 \nonumber \\
  + 3(7\varphi_{y}\psi_{y}r
 - 10\Delta)\varphi_{yyyy}\varphi_{y}^2r\Big],
 \end{eqnarray}
\begin{eqnarray}
 \fl  K_{4} =  - (\varphi_{y}^3\Delta)^{-1}\Big[(2(45r_{yy}r_{y}\Delta - 10r_{yy}\Delta_{x}
- 55r_{yy}\Delta_{y}r + 50r_{y}^2\Delta_{y} \nonumber \\
 - 20r_{y}\Delta_{xy} - 50r_{y}\Delta_{yy}r + 11\Delta_{xyy}r +
2\Delta_{xxy} +
17\Delta_{yyy}r^2 \nonumber \\
  - 20r_{yyy}r\Delta - 5r_{x}\Delta_{yy} -
10r_{xy}\Delta_{y} - 5r_{xyy}\Delta) \nonumber \\
  - 5(7\varphi_{y}^4\beta\psi r
+ 7\varphi_{y}^3\psi_{y}\alpha r - 3\varphi_{y}^2\alpha\Delta +
7\psi_{yyyy}r)\varphi_{y}r^2)\varphi_{y}^3 \nonumber \\
  + 15((3((5(\Delta_{x} +
5\Delta_{y}r) - 14r_{y}\Delta)r
 - r_{x}\Delta) - 35\varphi_{y}\psi_{yy}r^3)\varphi_{y}
 \nonumber \\
  + 35(\varphi_{y}\psi_{y}r - 3\Delta)\varphi_{yy}r^2)\varphi_{yy}^2
 - 10(\Delta_{xx} + 31\Delta_{yy}r^2 + 13\Delta_{xy}r
 \nonumber \\
  - 8(\Delta_{x} + 6\Delta_{y}r - 2r_{y}\Delta)r_{y}
 - 26r_{yy}r\Delta - 4r_{x}\Delta_{y} - 4r_{xy}\Delta
 \nonumber \\
  - 21\varphi_{y}\psi_{yyy}r^3)\varphi_{yy}\varphi_{y}^2 - 10(((5(\Delta_{x}
 + 5\Delta_{y}r) - 14r_{y}\Delta)r - r_{x}\Delta
 \nonumber \\
  - 14\varphi_{y}\psi_{yy}r^3)\varphi_{y} + 5(7\varphi_{y}\psi_{y}r
 - 18\Delta)\varphi_{yy}r^2)\varphi_{yyy}\varphi_{y}
 \nonumber \\
  + 5(7\varphi_{y}\psi_{y}r - 15\Delta)\varphi_{yyyy}\varphi_{y}^2r^2\Big], \label{eq:38}
\end{eqnarray}
\begin{eqnarray}
\fl K_{3} =  - (\varphi_{y}^3\Delta)^{-1}\Big[((13\Delta_{xxy}
+ 35\Delta_{yyy}r^2)r + \Delta_{xxx} + 31\Delta_{xyy}r^2\nonumber \\
  -
5(3\Delta_{xx} + 26\Delta_{yy}r^2 + 23\Delta_{xy}r - (15\Delta_{x}
+ 49\Delta_{y}r - 25r_{y}\Delta)r_{y})r_{y} \nonumber \\
  - 5(13\Delta_{x} +
32\Delta_{y}r - 50r_{y}\Delta)r_{yy}r - 65r_{yyy}r^2\Delta -
5(3\Delta_{xy} + 5\Delta_{yy}r \nonumber \\
  - 16r_{y}\Delta_{y} -
7r_{yy}\Delta)r_{x} - 5r_{xx}\Delta_{y} - 5r_{xxy}\Delta
 \nonumber \\
  - 5(3\Delta_{x} + 11\Delta_{y}r - 15r_{y}\Delta)r_{xy}
 - 30r_{xyy}r\Delta - 5(7\varphi_{y}^4\beta\psi r
 + 7\varphi_{y}^3\psi_{y}\alpha r \nonumber \\
  - 4\varphi_{y}^2\alpha\Delta
 + 7\psi_{yyyy}r)\varphi_{y}r^3)\varphi_{y}^3 - 5(2((2(2\Delta_{xx}
 + 17\Delta_{yy}r^2 + 11\Delta_{xy}r) \nonumber \\
  - (29\Delta_{x} + 75\Delta_{y}r
 - 51r_{y}\Delta)r_{y} - 45r_{yy}r\Delta)r - (3\Delta_{x}
 + 13\Delta_{y}r - 13r_{y}\Delta)r_{x} \nonumber \\
  - r_{xx}\Delta
 - 14r_{xy}r\Delta - 21\varphi_{y}\psi_{yyy}r^4)\varphi_{y}^2
 - 3((6((5(\Delta_{x} + 3\Delta_{y}r) - 13r_{y}\Delta)r
 \nonumber \\
  - 2r_{x}\Delta) - 35\varphi_{y}\psi_{yy}r^3)\varphi_{y}
 + 35(\varphi_{y}\psi_{y}r - 4\Delta)\varphi_{yy}r^2)\varphi_{yy}r)\varphi_{yy}
 \nonumber \\
  - 10(2((5(\Delta_{x} + 3\Delta_{y}r) - 13r_{y}\Delta)r - 2r_{x}\Delta
 - 7\varphi_{y}\psi_{yy}r^3)\varphi_{y}
 \nonumber \\
 + 5(7\varphi_{y}\psi_{y}r -
24\Delta)\varphi_{yy}r^2)\varphi_{yyy}\varphi_{y}r
 + 5(7\varphi_{y}\psi_{y}r
 - 20\Delta)\varphi_{yyyy}\varphi_{y}^2r^3\Big], \label{eq:39}
\end{eqnarray}
\begin{eqnarray}
\fl  K_{2} =   - (\varphi_{y}^3\Delta)^{-1}\Big[((3((5\Delta_{xxy} +
7\Delta_{yyy}r^2)r + \Delta_{xxx}
+ 7\Delta_{xyy}r^2) \nonumber \\
  - (3(13\Delta_{xx} + 28\Delta_{yy}r^2 +
39\Delta_{xy}r) + (204r_{y}\Delta - 161\Delta_{x} -
217\Delta_{y}r)r_{y})r_{y} \nonumber \\
  - (79\Delta_{x} + 116\Delta_{y}r -
264r_{y}\Delta)r_{yy}r - 54r_{yyy}r^2\Delta)r \nonumber \\
  - (3(2\Delta_{xx} +
7\Delta_{yy}r^2 + 11\Delta_{xy}r) + (171r_{y}\Delta - 64\Delta_{x}
- 140\Delta_{y}r)r_{y} \nonumber \\
  - 72r_{yy}r\Delta - 18r_{x}\Delta_{y})r_{x}
- (4\Delta_{x} + 11\Delta_{y}r - 21r_{y}\Delta)r_{xx} -
12r_{xxy}r\Delta \nonumber \\
  - r_{xxx}\Delta - ((37\Delta_{x} + 53\Delta_{y}r
- 150r_{y}\Delta)r - 33r_{x}\Delta)r_{xy} - 33r_{xyy}r^2\Delta \nonumber \\
  -
3(7\varphi_{y}^4\beta\psi r + 7\varphi_{y}^3\psi_{y}\alpha r -
5\varphi_{y}^2\alpha\Delta +
7\psi_{yyyy}r)\varphi_{y}r^4)\varphi_{y}^3 \nonumber \\
  - 3(2(5((2\Delta_{xx} +
7\Delta_{yy}r^2 + 6\Delta_{xy}r - (13\Delta_{x} + 19\Delta_{y}r -
20r_{y}\Delta)r_{y} \nonumber \\
  - 13r_{yy}r\Delta)r^2 - ((3\Delta_{x} +
5\Delta_{y}r - 11r_{y}\Delta)r - r_{x}\Delta)r_{x} - r_{xx}r\Delta
\nonumber \\
  - 6r_{xy}r^2\Delta) - 21\varphi_{y}\psi_{yyy}r^5)\varphi_{y}^2 -
15((2((5(\Delta_{x} + 2\Delta_{y}r) - 12r_{y}\Delta)r \nonumber \\
  -
3r_{x}\Delta) - 7\varphi_{y}\psi_{yy}r^3)\varphi_{y} +
7(\varphi_{y}\psi_{y}r -
5\Delta)\varphi_{yy}r^2)\varphi_{yy}r^2)\varphi_{yy} \nonumber \\
  -
2(2(5((5(\Delta_{x} + 2\Delta_{y}r) - 12r_{y}\Delta)r -
3r_{x}\Delta) - 21\varphi_{y}\psi_{yy}r^3)\varphi_{y} \nonumber \\
  +
15(7\varphi_{y}\psi_{y}r -
30\Delta)\varphi_{yy}r^2)\varphi_{yyy}\varphi_{y}r^2 +
3(7\varphi_{y}\psi_{y}r -
25\Delta)\varphi_{yyyy}\varphi_{y}^2r^4\Big], \label{eq:40}
\end{eqnarray}
\begin{eqnarray}
\fl K_{1}=  - (\varphi_{y}^3\Delta)^{-1}\Big[((7(\Delta_{xxy} +
\Delta_{yyy}r^2)r + 3\Delta_{xxx} + 7\Delta_{xyy}r^2 \nonumber \\
  -
(33\Delta_{xx}  + 28\Delta_{yy}r^2 + 49\Delta_{xy}r +
2(59r_{y}\Delta - 56\Delta_{x} - 42\Delta_{y}r)r_{y})r_{y} \nonumber \\
  -
(43\Delta_{x} + 42\Delta_{y}r - 128r_{y}\Delta)r_{yy}r -
23r_{yyy}r^2\Delta)r^2 \nonumber \\
  - ((12\Delta_{xx} + 7\Delta_{yy}r^2 +
21\Delta_{xy}r + 2(86r_{y}\Delta - 49\Delta_{x} -
35\Delta_{y}r)r_{y} \nonumber \\
  - 49r_{yy}r\Delta)r + (85r_{y}\Delta -
15\Delta_{x} - 21\Delta_{y}r)r_{x})r_{x} \nonumber \\
  - ((8\Delta_{x} +
7\Delta_{y}r - 32r_{y}\Delta)r - 10r_{x}\Delta)r_{xx} -
9r_{xxy}r^2\Delta - 2r_{xxx}r\Delta \nonumber \\
  - ((29\Delta_{x} +
21\Delta_{y}r - 95r_{y}\Delta)r - 46r_{x}\Delta)r_{xy}r -
16r_{xyy}r^3\Delta \nonumber \\
  - (7\varphi_{y}^4\beta\psi r +
7\varphi_{y}^3\psi_{y}\alpha r - 6\varphi_{y}^2\alpha\Delta +
7\psi_{yyyy}r)\varphi_{y}r^5)\varphi_{y}^3 \nonumber \\
  - (2(5((4\Delta_{xx} +
7\Delta_{yy}r^2 + 7\Delta_{xy}r - (23\Delta_{x} + 21\Delta_{y}r -
31r_{y}\Delta)r_{y} \nonumber \\
  - 17r_{yy}r\Delta)r^2 - ((9\Delta_{x} +
7\Delta_{y}r - 27r_{y}\Delta)r - 6r_{x}\Delta)r_{x} -
3r_{xx}r\Delta \nonumber \\
  - 10r_{xy}r^2\Delta) -
21\varphi_{y}\psi_{yyy}r^5)\varphi_{y}^2 - 15((3((5\Delta_{x} +
7\Delta_{y}r - 11r_{y}\Delta)r \nonumber \\
  - 4r_{x}\Delta) -
7\varphi_{y}\psi_{yy}r^3)\varphi_{y} + 7(\varphi_{y}\psi_{y}r -
6\Delta)\varphi_{yy}r^2)\varphi_{yy}r^2)\varphi_{yy}r \nonumber \\
  -
2((5((5\Delta_{x} + 7\Delta_{y}r - 11r_{y}\Delta)r - 4r_{x}\Delta)
- 14\varphi_{y}\psi_{yy}r^3)\varphi_{y} \nonumber \\
  + 5(7\varphi_{y}\psi_{y}r
- 36\Delta)\varphi_{yy}r^2)\varphi_{yyy}\varphi_{y}r^3 +
(7\varphi_{y}\psi_{y}r -
30\Delta)\varphi_{yyyy}\varphi_{y}^2r^5\Big], \label{eq:41}
\end{eqnarray}
\begin{eqnarray}
\fl K_{0} =  (\varphi_{y}^3\Delta)^{-1}\Big[((((2(r_{xxy} +
2r_{yyy}r^2)r + r_{xxx} + 3r_{xyy}r^2)\Delta \nonumber \\
  + 3(3\Delta_{x} +
2\Delta_{y}r - 8r_{y}\Delta)r_{yy}r^2)r - ((10r_{x} +
11r_{y}r)\Delta \nonumber \\
  - (4\Delta_{x} + \Delta_{y}r)r)r_{xx} - ((13r_{x}
+ 20r_{y}r)\Delta
 - (7\Delta_{x} + 3\Delta_{y}r)r)r_{xy}r \nonumber \\
  + ((\varphi_{y}^4\beta\psi
 + \varphi_{y}^3\psi_{y}\alpha + \psi_{yyyy})r
 - \varphi_{y}^2\alpha\Delta)\varphi_{y}r^5 + (9\Delta_{xx}
 + 4\Delta_{yy}r^2 \nonumber \\
  + 7\Delta_{xy}r - 2(13\Delta_{x} + 6\Delta_{y}r
 - 12r_{y}\Delta)r_{y})r_{y}r^2 - ((\Delta_{xxy} + \Delta_{yyy}r^2)r
 \nonumber \\
  + \Delta_{xxx} + \Delta_{xyy}r^2)r^2)r
 - ((2((17\Delta_{x} + 5\Delta_{y}r - 23r_{y}\Delta)r_{y}
 + 6r_{yy}r\Delta) \nonumber \\
  - (6\Delta_{xx} + \Delta_{yy}r^2
 + 3\Delta_{xy}r))r^2 - (5(3r_{x} + 8r_{y}r)\Delta
 \nonumber \\
  - 3(5\Delta_{x} + \Delta_{y}r)r)r_{x})r_{x})\varphi_{y}^3
 - ((2((5(r_{xx} + 3r_{yy}r^2 + 2r_{xy}r)\Delta
 \nonumber \\
  + 3\varphi_{y}\psi_{yyy}r^4 + 5(5\Delta_{x}
 + 3\Delta_{y}r - 6r_{y}\Delta)r_{y}r - 5(\Delta_{xx}
 + \Delta_{yy}r^2 + \Delta_{xy}r)r)r \nonumber \\
  - 5((3r_{x}
 + 7r_{y}r)\Delta - (3\Delta_{x} + \Delta_{y}r)r)r_{x})\varphi_{y}^2
 - 15((3(r_{x} + 2r_{y}r)\Delta \nonumber \\
  + \varphi_{y}\psi_{yy}r^3
 - 3(\Delta_{x} + \Delta_{y}r)r)\varphi_{y}
 - (\varphi_{y}\psi_{y}r - 7\Delta)\varphi_{yy}r^2)\varphi_{yy}r^2)\varphi_{yy}
 \nonumber \\
  + (2((5(r_{x} + 2r_{y}r)\Delta + 2\varphi_{y}\psi_{yy}r^3
 - 5(\Delta_{x} + \Delta_{y}r)r)\varphi_{y} \nonumber \\
  - 5(\varphi_{y}\psi_{y}r
 - 6\Delta)\varphi_{yy}r^2)\varphi_{yyy} + (\varphi_{y}\psi_{y}r
 - 5\Delta)\varphi_{yyyy}\varphi_{y}r^2)\varphi_{y}r^2)r^2\Big]. \label{eq:42}
\end{eqnarray}
{\bf Definition 2}.
 We call  (\ref{eq:21}) with arbitrary coefficients
 $r = r(x, y)$,\linebreak  $F_0 = F_0(x, y), ~F_1 = F_1(x, y), ~F_2 = F_2(x, y),
  ~H_0 = H_0(x, y), ~H_1 = H_1(x, y),$\linebreak $ ~H_2 = H_2(x, y), ~J_0 = J_0(x, y),
 ~J_1 = J_1(x, y), ~J_2 = J_2(x, y),  ~J_3 = J_3(x, y), ~J_4 = J_4(x, y),$  and
 $K_i = K_i(x, y), ~(i = 0, \ldots, 7),$ the
second candidate for linearization.

Thus, we showed that every linearizable fourth-order
 equations belong either to the class of  (\ref{eq:09})
 or to the class of  (\ref{eq:21}).
 In Sections \ref{candidate2} and  \ref{candidate1}, we formulate
 the main theorems containing necessary and sufficient conditions for linearization
as well as  the methods for constructing the linearizing point
transformations for each candidate. Proofs of the main theorems and
illustrative examples are provided in the subsequent sections.

\subsection{The linearization test for equation (\ref{eq:09})}
 \label{candidate2}

Consider the first canditate for linearization, i.e. equation
(\ref{eq:09}). In this case, the linearizing transformations
(\ref{eq:02}) have the form
\begin{equation}
 \label{eq:43}
 t=\varphi (x),\quad u=\psi (x,y).\\[1ex]
 \end{equation}
{\bf Theorem 1}.
 \label{eq:09th}
  Equation (\ref{eq:09})
$$
  \begin{array}{ll}
  y^{(4)}  & + (A_{1}y' + A_{0})y''' + B_{0}y''^2 + (C_{2}y'^2 + C_{1}y' + C_{0})y''
 \\ &+ D_{4}y'^4 + D_{3}y'^3 + D_{2}y'^2 + D_{1}y' + D_{0} = 0,
 \end{array}
$$
is linearizable if and only if its coefficients obey the following
ten equations
\begin{eqnarray}
 \fl A_{0y} - A_{1x} = 0,\label{eq:45}\\[1.5ex]
 \fl 4B_{0} - 3A_{1} = 0, \label{eq:46}\\[1.5ex]
 \fl 12A_{1y} + 3A_{1}^2 - 8C_{2} = 0, \label{eq:47}\\[1.5ex]
 \fl 12A_{1x} + 3A_{0}A_{1} - 4C_{1} = 0, \label{eq:48}\\[1.5ex]
 \fl 32C_{0y} + 12A_{0x}A_{1} - 16C_{1x} + 3A_{0}^2A_{1} -
4A_{0}C_{1}
= 0, \label{eq:49}\\[1.5ex]
 \fl 4C_{2y} + A_{1}C_{2} - 24D_{4} = 0, \label{eq:50}\\[1.5ex]
 \fl 4C_{1y} + A_{1}C_{1} - 12D_{3} = 0, \label{eq:51}\\[1.5ex]
 \fl 16C_{1x} - 12A_{0x}A_{1} - 3A_{0}^2A_{1} + 4A_{0}C_{1} +
8A_{1}C_{0} - 32D_{2}=0, \label{eq:52}\\[1.5ex]
\nonumber
  \fl 192D_{2x} + 36A_{0x}A_{0}A_{1} - 48A_{0x}C_{1} -
48C_{0x}A_{1} - 288D_{1y} + 9A_{0}^3A_{1}   - 12A_{0}^2C_{1}
\nonumber \\- 36A_{0}A_{1}C_{0} + 48A_{0}D_{2} + 32C_{0}C_{1} =
0, \label{eq:53}\\[1.5ex]
\nonumber
 \fl 384D_{1xy} - \Big[3((3A_{0}A_{1} - 4C_{1})A_{0}^2 +
16(2A_{1}D_{1} + C_{0}C_{1})
 - 16(A_{1}C_{0} - D_{2})A_{0})A_{0} \nonumber \\- 32(4(C_{1}D_{1} -
2C_{2}D_{0} + C_{0}D_{2})
 + (3A_{1}D_{0} - C_{0}^2)A_{1}) - 96D_{1y}A_{0} \nonumber \\+ 384D_{0y}A_{1}
+ 1536D_{0yy}
 - 16(3A_{0}A_{1} - 4C_{1})C_{0x} \nonumber \\+ 12((3A_{0}A_{1} -
4C_{1})A_{0}
 - 4(A_{1}C_{0}  - 4D_{2}))A_{0x}\Big] = 0. \label{eq:54}
\end{eqnarray}
Provided that the conditions (\ref{eq:45})-(\ref{eq:54}) are
satisfied, the linearizing transformation (\ref{eq:43}) is defined
by a fourth-order ordinary differential equation for the function
$\varphi (x),$ namely by the Riccati equation
 \begin{equation}
 \label{eq:55}
 40\frac{d \chi}{dx} - 20\chi^2 = 8C_0 -3A_0^2 -12A_{0x},
 \end{equation}
 for
 \begin{equation}
 \label{eq:56}
 \chi = \frac{\varphi_{xx}}{\varphi_x}\,,
 \end{equation}
and by the following integrable system of partial differential
equations for $\psi (x,y)$
\begin{equation}\label{eq:57}
    4\psi_{yy} =  \psi_{y}A_{1},
\end{equation}
\begin{equation}\label{eq:58}
    4\psi_{xy} =  \psi_{y}(A_{0} + 6\chi),
\end{equation}
and
\begin{eqnarray}
 \nonumber
\fl 1600\psi_{xxxx} =   9600\psi_{xxx}\chi + 160\psi_{xx}( -
12A_{0x} - 3A_{0}^2 - 90\chi^2 + 8C_{0}) \nonumber \\  +
40\psi_{x}(12A_{0x}A_{0} + 72A_{0x}\chi - 16C_{0x} + 3A_{0}^3 +
18A_{0}^2\chi  - 12A_{0}C_{0} \nonumber \\  + 120\chi^3 - 48\chi
C_{0} + 24D_{1} - 8\Omega)  + \psi(144A_{0x}^2 + 72A_{0x}A_{0}^2 -
352A_{0x}C_{0} \nonumber \\ - 160C_{0xx} - 80C_{0x}A_{0}
 - 1600D_{0y} + 640D_{1x} - 80\Omega_{x} + 9A_{0}^4
- 88A_{0}^2C_{0} \nonumber \\  + 160A_{0}D_{1}  + 30A_{0}\Omega -
400A_{1}D_{0} + 300\chi\Omega + 144C_{0}^2) + 1600\psi_{y}D_{0},
\label{eq:59}
\end{eqnarray}
where $\chi$ is given by (\ref{eq:56}) and $\Omega$ is the
 following expression
\begin{equation}
 \label{eq:60}
 \Omega = A_{0}^3 - 4A_{0}C_{0} + 8D_{1} - 8C_{0x}
+ 6A_{0x}A_{0} + 4A_{0xx}.
\end{equation}
Finally, the coefficients $\alpha $ and $\beta$ of the resulting
linear equation (\ref{eq:03}) is given by
\begin{equation}\label{eq:61}
    \alpha =  \frac{\Omega}{8\varphi_{x}^3} ~,
\end{equation}
and
\begin{eqnarray}
\fl \beta =  (1600\varphi_{x}^4)^{-1}( - 144A_{0x}^2 -
72A_{0x}A_{0}^2 + 352A_{0x}C_{0} + 160C_{0xx} + 80C_{0x}A_{0}
\nonumber \\  + 1600D_{0y} - 640D_{1x} + 80\Omega_{x} - 9A_{0}^4 +
88A_{0}^2C_{0} - 160A_{0}D_{1} - 30A_{0}\Omega \nonumber \\ +
400A_{1}D_{0} - 300\chi\Omega - 144C_{0}^2).\label{eq:62}
\end{eqnarray}
{\bf Remark 1}.
 Since the system of equations (\ref{eq:45})-(\ref{eq:54}) provides
 the necessary and sufficient conditions for linearization,
 it is invariant with respect to the transformations (\ref{eq:43}).
 It means that the left-hand sides
 of  (\ref{eq:45})-(\ref{eq:54}) are {\it relative
 invariants} (of the second-order) for the equivalence group (\ref{eq:43}).

\subsection{The linearization test for equation (\ref{eq:21})}
\label{candidate1} The following theorem provides the test for
linearization of the
 second candidate. The necessary and
sufficient conditions comprise {\it eighteen} differential equations
(\ref{eq:64})-(\ref{eq:82}) for {\it twenty one} coefficients of the
 (\ref{eq:21}). The linearizing change of variables
(\ref{eq:02}) is determined by  (\ref{eq:83})-(\ref{eq:86})
for the functions $\varphi(x, y)$ and $\psi(x, y).$\\[1ex]
{\bf Theorem 2}.
 \label{eq:21th}
  Equation (\ref{eq:21})
$$
\begin{array}{c}
  y^{\left( 4 \right)}   +  \frac{1}
{{y'  +  r}}(  - 10 y''  +  F_2 y'^2  +  F_1 y'  +  F_{0 } )y'''  \hfill \\
   + \frac{1}
{{( y'  +  r)^2 }}\left[
  15y''^3   +  ( H_2 y'^2  +  H_1 y' +  H_0 )y''^2  \hfill \right.
\\
   +  (J_4 y'^4  +  J_3 y'^3  +  J_2 y'^2  +  J_1 y'  + J_{0 } )y'' \hfill \\
   +  K_7 y'^7  +  K_6 y'^6  +  K_5 y'^5  +  K_4 y'^4  \hfill \\
\left.   +  K_3 y'^3  +  K_2 y'^2  +  K_1 y' +  K_0  \hfill
\right] = 0, \hfill \\
\end{array}
$$
  is linearizable if and only if its coefficients
 obey the following equations
 \begin{eqnarray}
& \fl 10r_{yy} =  - (F_{1y} + F_{2x} + F_{2y}r + r_{y}F_{2}),\label{eq:64}\\[1.5ex]
& \fl 10r_{x} =  10r_{y}r - F_{0} + F_{1}r - F_{2}r^2,\label{eq:65}\\[1.5ex]
& \fl H_{2} =  - 3F_{2}, \label{eq:67}\\[1.5ex]
& \fl 4H_{1} =  - 3(5F_{1} - 2F_{2}r), \label{eq:68}\\[1.5ex]
& \fl 4H_{0} =  - 3(6F_{0} - F_{1}r), \label{eq:69}\\[1.5ex]
 & \fl 10F_{1yy} =  - (F_{1y}F_{2} - 40F_{2xy}
- 16F_{2x}F_{2} + 20F_{2yy}r + 40F_{2y}r_{y}  + 14F_{2y}F_{2}r
 + 20J_{4x} \nonumber \\ &- 20J_{4y}r
+ 14r_{y}F_{2}^2 - 40r_{y}J_{4}), \label{eq:70}\\[1.5ex]
& \fl 12F_{2x} =  12F_{2y}r - 3F_{1}F_{2} + 6F_{2}^2r + 4J_{3} -
16J_{4}r, \label{eq:71}\\[1.5ex]
& \fl 60F_{1x} =  60F_{1y}r - 36F_{0}F_{2} - 15F_{1}^2 +
66F_{1}F_{2}r - 36F_{2}^2r^2 + 40J_{2}  - 80J_{3}r  + 80J_{4}r^2,
\label{eq:72}\\[1.5ex]
 & \fl 60F_{0x} =  60F_{0y}r - 51F_{0}F_{1} + 66F_{0}F_{2}r +
36F_{1}^2r - 72F_{1}F_{2}r^2 + 36F_{2}^2r^3  +
60J_{1} \nonumber \\ &- 80J_{2}r + 80J_{3}r^2 - 80J_{4}r^3, \label{eq:73}\\[1.5ex]
 & \fl 20J_{0} =  9F_{0}^2 - 18F_{0}F_{1}r + 18F_{0}F_{2}r^2
+ 9F_{1}^2r^2 - 18F_{1}F_{2}r^3 + 9F_{2}^2r^4  +
20J_{1}r \nonumber \\ & - 20J_{2}r^2 + 20J_{3}r^3 - 20J_{4}r^4, \label{eq:74}\\[1.5ex]
& \fl 120J_{3yy} =  216F_{1y}F_{2y} + 54F_{1y}F_{2}^2 -
48F_{1y}J_{4} + 360F_{2yy}r_{y} + 90F_{2yy}F_{1}  - 180F_{2yy}F_{2}r
\nonumber \\ &- 432F_{2y}^2r + 324F_{2y}r_{y}F_{2} +
189F_{2y}F_{1}F_{2} - 486F_{2y}F_{2}^2r  - 192F_{2y}J_{3} \nonumber
\\ &+ 864F_{2y}J_{4}r - 60J_{3y}F_{2} + 720J_{4xy} + 180J_{4x}F_{2}
- 240J_{4yy}r  \nonumber \\ &- 1200J_{4y}r_{y} + 60J_{4y}F_{2}r +
720K_{6x} - 720K_{6y}r - 5040K_{7x}r \nonumber \\ &+ 5040K_{7y}r^2 +
36r_{y}F_{2}^3 - 432r_{y}F_{2}J_{4} - 2160r_{y}K_{6} +
15120r_{y}K_{7}r \nonumber \\ &+ 504F_{0}K_{7} + 36F_{1}F_{2}^3  -
102F_{1}F_{2}J_{4} - 504F_{1}K_{7}r - 72F_{2}^4r \nonumber \\ &-
48F_{2}^2J_{3} + 396F_{2}^2J_{4}r + 504F_{2}K_{7}r^2 + 136J_{3}J_{4}
- 544J_{4}^2r, \label{eq:75}
\end{eqnarray}
\begin{eqnarray}
\nonumber
 & \fl 240J_{4xyy} =  - (36F_{1y}F_{2yy} +
162F_{1y}F_{2y}F_{2} - 72F_{1y}J_{4y} + 36F_{1y}F_{2}^3 -
168F_{1y}F_{2}J_{4} \nonumber
\\ & - 72F_{1y}K_{6}  - 168F_{1y}K_{7}r -
72F_{2yy}F_{2y}r + 144F_{2yy}r_{y}F_{2} \nonumber \\ & +
54F_{2yy}F_{1}F_{2} - 108F_{2yy}F_{2}^2r - 72F_{2yy}J_{3} +
288F_{2yy}J_{4}r + 432F_{2y}^2r_{y} \nonumber \\ & +
108F_{2y}^2F_{1} - 540F_{2y}^2F_{2}r - 144F_{2y}J_{3y} +
528F_{2y}J_{4x} + 192F_{2y}J_{4y}r \nonumber \\ & +
324F_{2y}r_{y}F_{2}^2 - 1008F_{2y}r_{y}J_{4} + 162F_{2y}F_{1}F_{2}^2
- 132F_{2y}F_{1}J_{4} \nonumber \\ &- 396F_{2y}F_{2}^3r  -
180F_{2y}F_{2}J_{3} + 1320F_{2y}F_{2}J_{4}r + 144F_{2y}K_{6}r
\nonumber \\ &- 336F_{2y}K_{7}r^2 - 36J_{3y}F_{2}^2  +
176J_{3y}J_{4} + 120J_{4xy}F_{2} + 132J_{4x}F_{2}^2 \nonumber \\
&- 432J_{4x}J_{4} - 240J_{4yyy}r - 960J_{4yy}r_{y}  -
120J_{4yy}F_{2}r - 768J_{4y}r_{y}F_{2} \nonumber \\ &-
138J_{4y}F_{1}F_{2} + 288J_{4y}F_{2}^2r + 184J_{4y}J_{3} -
1008J_{4y}J_{4}r + 960K_{6xy} \nonumber
\\ & + 240K_{6x}F_{2} - 960K_{6yy}r -
3840K_{6y}r_{y}  - 240K_{6y}F_{2}r - 1920K_{7xy}r \nonumber
\\ & - 2400K_{7xx} +
2880K_{7x}r_{y} - 600K_{7x}F_{1}  - 480K_{7x}F_{2}r + 4320K_{7yy}r^2
\nonumber \\ &+ 24000K_{7y}r_{y}r + 432K_{7y}F_{0} + 168K_{7y}F_{1}r
+ 912K_{7y}F_{2}r^2 \nonumber \\
& + 20160r_{y}^2K_{7} + 1728r_{y}F_{1}K_{7} + 36r_{y}F_{2}^4 -
264r_{y}F_{2}^2J_{4}  - 1248r_{y}F_{2}K_{6} \nonumber \\ &+
5280r_{y}F_{2}K_{7}r + 160r_{y}J_{4}^2 + 408F_{0}F_{2}K_{7} +
150F_{1}^2K_{7} + 27F_{1}F_{2}^4 \nonumber \\
& - 120F_{1}F_{2}^2J_{4} - 168F_{1}F_{2}K_{6} + 168F_{1}F_{2}K_{7}r
- 54F_{2}^5r - 36F_{2}^3J_{3} \nonumber \\ & + 384F_{2}^3J_{4}r +
336F_{2}^2K_{6}r - 1344F_{2}^2K_{7}r^2 + 160F_{2}J_{3}J_{4} -
640F_{2}J_{4}^2r \nonumber \\ & - 400J_{2}K_{7} + 224J_{3}K_{6} -
368J_{3}K_{7}r - 896J_{4}K_{6}r + 3872J_{4}K_{7}r^2 \nonumber \\ & +
672F_{0y}K_{7}), \label{eq:76}
\end{eqnarray}
\begin{eqnarray}
& \fl 4J_{4x} =  4J_{4y}r - F_{1}J_{4} + 2F_{2}J_{4}r - 4K_{5} +
24K_{6}r - 84K_{7}r^2, \label{eq:77}
\end{eqnarray}
\begin{eqnarray}
\nonumber &  \fl 60F_{0yy} =  - (30F_{0y}F_{2} + 36F_{1y}F_{1} -
36F_{1y}F_{2}r - 60F_{2yy}r^2 + 24F_{2y}F_{0}  - 36F_{2y}F_{1}r
\nonumber \\ &- 54F_{2y}F_{2}r^2 - 40J_{2y} + 40J_{3y}r +
80J_{4y}r^2  - 36r_{y}F_{1}F_{2} + 36r_{y}F_{2}^2r \nonumber \\ &+
40r_{y}J_{3} - 80r_{y}J_{4}r + 6F_{0}F_{2}^2 - 6F_{0}J_{4} +
9F_{1}^2F_{2} - 18F_{1}F_{2}^2r \nonumber
\\ & - 12F_{1}J_{3} + 24F_{1}J_{4}r -
6F_{2}^3r^2 - 10F_{2}J_{2}  + 22F_{2}J_{3}r + 26F_{2}J_{4}r^2
\nonumber \\ &- 60K_{4} + 180K_{5}r - 180K_{6}r^2 - 420K_{7}r^3),
\label{eq:78}
\end{eqnarray}
\begin{eqnarray}
\nonumber
 &  \fl 20J_{2x} =  20J_{2y}r + 20J_{3x}r - 20J_{3y}r^2 -
14F_{0}J_{3} + 28F_{0}J_{4}r - 5F_{1}J_{2}  + 19F_{1}J_{3}r
\nonumber \\ &- 28F_{1}J_{4}r^2 + 10F_{2}J_{2}r - 24F_{2}J_{3}r^2 +
28F_{2}J_{4}r^3  - 120K_{3} + 360K_{4}r \nonumber \\ &- 640K_{5}r^2
+ 840K_{6}r^3 - 840K_{7}r^4, \label{eq:79}
\end{eqnarray}
\begin{eqnarray}
\nonumber
 &  \fl 60J_{1x} =  60J_{1y}r - 40J_{3x}r^2 + 40J_{3y}r^3 -
42F_{0}J_{2} + 42F_{0}J_{3}r - 70F_{0}J_{4}r^2  - 15F_{1}J_{1}
\nonumber \\ &+ 42F_{1}J_{2}r - 52F_{1}J_{3}r^2 + 70F_{1}J_{4}r^3 +
30F_{2}J_{1}r  - 42F_{2}J_{2}r^2 \nonumber \\ &+ 62F_{2}J_{3}r^3 -
70F_{2}J_{4}r^4 - 600K_{2} + 1080K_{3}r - 1380K_{4}r^2 \nonumber
\\ & + 1700K_{5}r^3 - 2100K_{6}r^4 + 2100K_{7}r^5, \label{eq:80}
\end{eqnarray}
\begin{eqnarray}
\nonumber & \fl 80K_{1} =  3F_{0}^2F_{1} - 6F_{0}^2F_{2}r -
6F_{0}F_{1}^2r + 18F_{0}F_{1}F_{2}r^2 - 12F_{0}F_{2}^2r^3 -
8F_{0}J_{1} \nonumber
\\ & + 16F_{0}J_{2}r - 24F_{0}J_{3}r^2 + 32F_{0}J_{4}r^3 +
3F_{1}^3r^2 - 12F_{1}^2F_{2}r^3 + 15F_{1}F_{2}^2r^4 \nonumber \\
& + 8F_{1}J_{1}r - 16F_{1}J_{2}r^2 + 24F_{1}J_{3}r^3 -
32F_{1}J_{4}r^4 - 6F_{2}^3r^5 - 8F_{2}J_{1}r^2 \nonumber \\ & +
16F_{2}J_{2}r^3 - 24F_{2}J_{3}r^4 + 32F_{2}J_{4}r^5 + 160K_{2}r -
240K_{3}r^2 + 320K_{4}r^3 \nonumber \\ & - 400K_{5}r^4 + 480K_{6}r^5
- 560K_{7}r^6, \label{eq:81}
\end{eqnarray}
\begin{eqnarray}
\nonumber
 & \fl 400K_{0} =  - (6F_{0}^3 - 33F_{0}^2F_{1}r +
48F_{0}^2F_{2}r^2 + 48F_{0}F_{1}^2r^2 - 126F_{0}F_{1}F_{2}r^3
 + 78F_{0}F_{2}^2r^4 \nonumber \\ &+ 40F_{0}J_{1}r - 80F_{0}J_{2}r^2
+ 120F_{0}J_{3}r^3 - 160F_{0}J_{4}r^4 - 21F_{1}^3r^3 \nonumber \\ &+
78F_{1}^2F_{2}r^4 - 93F_{1}F_{2}^2r^5 - 40F_{1}J_{1}r^2 +
80F_{1}J_{2}r^3 - 120F_{1}J_{3}r^4 \nonumber \\ &+ 160F_{1}J_{4}r^5
+ 36F_{2}^3r^6 + 40F_{2}J_{1}r^3 - 80F_{2}J_{2}r^4 +
120F_{2}J_{3}r^5 \nonumber \\ &- 160F_{2}J_{4}r^6 - 400K_{2}r^2  +
800K_{3}r^3 - 1200K_{4}r^4 + 1600K_{5}r^5 \nonumber
\\ & - 2000K_{6}r^6 +
2400K_{7}r^7). \label{eq:82}
\end{eqnarray}
Provided that the conditions (\ref{eq:64})-(\ref{eq:82}) are
satisfied, the transformations (\ref{eq:02}) mapping equation
(\ref{eq:21}) to a linear equation (\ref{eq:03}) is obtained by
solving the following compatible system of equations for the
functions $\varphi (x,y)$ and $\psi (x,y)$
\begin{equation}\label{eq:83}
    \varphi_{x} = r\varphi_{y},
\end{equation}
\begin{equation}\label{eq:84}
    \varphi_{y}\psi_{x} = r\varphi_{y}\psi_{y}  - \Delta,
\end{equation}
\begin{equation}\label{eq:85}
    10\Delta \varphi_{yy} = \varphi_{y}(4\Delta_{y} - F_{2}\Delta),
\end{equation}
and
 \begin{eqnarray}
 \fl 500 \varphi_{y} \psi_{yyyy} \Delta^3 = &300 \psi_{yyy}
\varphi_{y} \Delta^2 (4 \Delta_{y} - F_{2}\Delta) +5 \psi_{yy}
\varphi_{y} \Delta ( - 120 F_{2y} \Delta^2  - 144 \Delta_{y}^2
\nonumber \\&+ 72 \Delta_{y} F_{2} \Delta - 39 F_{2}^2 \Delta^2 + 80
J_{4} \Delta^2) +\psi_{y} \varphi_{y} ( - 500 \varphi_{y}^3 \alpha
\Delta^3 \nonumber \\& -150 F_{2yy} \Delta^3 + 360 F_{2y} \Delta_{y}
\Delta^2 -165 F_{2y} F_{2} \Delta^3 + 100 J_{4y} \Delta^3
 + 96 \Delta_{y}^3 \nonumber \\& - 72 \Delta_{y}^2 F_{2} \Delta + 108 \Delta_{y}
F_{2}^2 \Delta^2 -240 \Delta_{y} J_{4} \Delta^2 - 24 F_{2}^3
\Delta^3  + 60  F_{2} J_{4} \Delta^3) \nonumber \\& -500 \psi
\varphi_{y}^5 \beta \Delta^3 +500 K_{7} \Delta^4.
 \label{eq:86}
 \end{eqnarray}
The coefficients $\alpha$ and $\beta$ of the resulting linear
equation (\ref{eq:03}) is given  by
\begin{equation}\label{eq:88}
    \alpha = \frac{\Theta}{8\varphi_{y}^3}  ~,
\end{equation}
and
 \begin{eqnarray}
 \fl \beta  = (1600\Delta\varphi_{y}^4)^{-1}\Big[\Delta( - 144F_{2y}^2 - 72F_{2y}F_{2}^2 +
352F_{2y}J_{4} + 160J_{4yy}
  + 80J_{4y}F_{2} \nonumber \\+ 640K_{6y} - 1600K_{7x} - 2880K_{7y}r +
80\Theta_{y} - 4480r_{y}K_{7}   - 400F_{1}K_{7} \nonumber \\-
9F_{2}^4 + 88F_{2}^2J_{4} + 160F_{2}K_{6} - 320F_{2}K_{7}r -
144J_{4}^2)
 - 120\Delta_{y}\Theta\Big],
 \label{eq:89}
 \end{eqnarray}
where  $\Theta$ is the following expression
\begin{equation}\label{eq:87}
    \Theta = (F_{2}^2 - 4J_{4})F_{2} - 8(K_{6} - 7K_{7} r)
- 8J_{4y} +6F_{2y}F_{2} + 4 F_{2yy}.
\end{equation}
{\bf Remark 2}.
 The equations (\ref{eq:64})-(\ref{eq:82}) define eighteen {\it relative
 invariants} of the third-order for the general
 point transformation group (\ref{eq:02}).

 \section{Proof of the linearization theorems}
The proof of the linearization theorems formulated above
 requires investigation of integrability conditions for the equations
 given in Section \ref{candidates}. We will consider the problem for
 the candidates (\ref{eq:09}) and (\ref{eq:21}) separately.
 The problem is formulated as follows. Given the coefficients $A_i(x,y),
 B_i(x,y),
 C_i(x,y), D_i(x,y)$ and  $F_i(x,y), H_i(x,y), J_i(x,y), K_i(x,y)$
  of the equations (\ref{eq:09}) and (\ref{eq:21}), respectively,
 find the integrability conditions of the respective equations for the functions
 $\varphi$ and $\psi.$

 \subsection{Proof of Theorem 1}
Let us turn to the proof of Theorem 1 on linearization of
(\ref{eq:09}). Namely, given the coefficients $A_i(x,y), B_i(x,y),
C_i(x,y), D_i(x,y)$ of  (\ref{eq:09}), we have to find the necessary
and sufficient conditions for integrability of the over-determined
system (\ref {eq:10})-(\ref{eq:20}) for the unknown functions
$\varphi (x)$ and $\psi (x,y).$

We first rewrite the expressions (\ref{eq:10}) and (\ref{eq:11}) for
$A_{1}$ and $A_{0}$ in the following form
\begin{equation}\label{pf1:01}
    \psi_{yy} = \frac{\psi_{y}A_{1}}{4} ,
\end{equation}
\begin{equation}\label{pf1:02}
    \psi_{xy} = \frac{(6\varphi_{xx} +
    \varphi_{x}A_{0})}{4\varphi_{x}}\psi_{y}.
\end{equation}
Comparing the mixed derivative $(\psi_{yy})_{x}=(\psi_{xy})_{y}$,
one arrives at  (\ref{eq:45})
    \[A_{0y} = A_{1x}.\]
Then (\ref{eq:12}), (\ref{eq:13}) and (\ref{eq:14}) are written in
the form
   \[ 3A_{1} - 4B_{0}=0,\]
    \[3A_{1}^2 - 8C_{2} + 12A_{1y}=0,\]
and \[    12A_{1x} + 3A_{0}A_{1} -
    4C_{1}=0,\]
respectively.  So that one obtains (\ref{eq:46}), (\ref{eq:47}) and
(\ref{eq:48}) respectively. Furthermore,  (\ref{eq:15}) for $C_{0}$
becomes
\begin{equation}\label{pf1:07}
    \varphi_{xxx} =  - \frac{(12A_{0x}\varphi_{x}^2 - 60\varphi_{xx}^2
+ 3\varphi_{x}^2A_{0}^2 - 8\varphi_{x}^2C_{0})}{40\varphi_{x}}.
\end{equation}
Differentiation of  (\ref{pf1:07}) with respect to $y$ yields
\[
    12A_{0x}A_{1} + 32C_{0y} - 16C_{1x} + 3A_{0}^2A_{1}
- 4A_{0}C_{1}=0.
\]
Thus one gets  (\ref{eq:49}). Therefore  (\ref{eq:16}),
(\ref{eq:17}) and (\ref{eq:18}) can be written in the form of
(\ref{eq:50}), (\ref{eq:51}) and (\ref{eq:52}), respectively.

One can determine $\alpha$ from  (\ref{eq:19}), as the following
\begin{equation}\label{pf1:12}
    \alpha = \frac{4A_{0xx} + 6A_{0x}A_{0} - 8C_{0x} + A_{0}^3 - 4A_{0}C_{0}
+ 8D_{1}}{8\varphi_{x}^3}.
\end{equation}
Since $\varphi=\varphi(x)$, we have $\alpha_{y}=0$ yields
(\ref{eq:53})
$$
 \begin{array}{ll}
    D_{2x} = & - \frac{1}{192}\Big[36A_{0x}A_{0}A_{1} - 48A_{0x}C_{1}
- 48C_{0x}A_{1} - 288D_{1y} + 9A_{0}^3A_{1}\\
 &- 12A_{0}^2C_{1} -
36A_{0}A_{1}C_{0} + 48A_{0}D_{2} + 32C_{0}C_{1}\Big].
\end{array}
$$
From  (\ref{eq:20}) one finds
\begin{eqnarray}
\psi_{xxxx}= & -
\frac{1}{40\varphi_{x}^3}\Big[32A_{0xx}\varphi_{x}^3\psi_{x} -
72A_{0x}\varphi_{xx}\varphi_{x}^2\psi_{x} +
48A_{0x}\varphi_{x}^3\psi_{xx}\nonumber \nonumber \\
& + 36A_{0x}\varphi_{x}^3\psi_{x}A_{0} -
48C_{0x}\varphi_{x}^3\psi_{x} - 120\varphi_{xx}^3\psi_{x} +
360\varphi_{xx}^2\varphi_{x}\psi_{xx} \nonumber \\ & -
240\varphi_{xx}\varphi_{x}^2\psi_{xxx} -
18\varphi_{xx}\varphi_{x}^2\psi_{x}A_{0}^2 +
48\varphi_{xx}\varphi_{x}^2\psi_{x}C_{0} + 40\varphi_{x}^7\beta\psi
\nonumber \\ & + 12\varphi_{x}^3\psi_{xx}A_{0}^2 -
32\varphi_{x}^3\psi_{xx}C_{0} + 5\varphi_{x}^3\psi_{x}A_{0}^3 -
20\varphi_{x}^3\psi_{x}A_{0}C_{0} \nonumber \\ & +
40\varphi_{x}^3\psi_{x}D_{1} - 40\varphi_{x}^3\psi_{y}D_{0}\Big].
 \label{pf1:14}
\end{eqnarray}
Forming the mixed derivative $(\psi_{xxxx})_{y}=(\psi_{xy})_{xxx}$
one obtains
 \begin{eqnarray}
\beta = & \frac{1}{1600\varphi_{x}^5}\Big[320A_{0xxx}\varphi_{x} -
1200A_{0xx}\varphi_{xx} + 360A_{0xx}\varphi_{x}A_{0} +
336A_{0x}^2\varphi_{x}\nonumber \\ & - 1800A_{0x}\varphi_{xx}A_{0} -
12A_{0x}\varphi_{x}A_{0}^2 + 32A_{0x}\varphi_{x}C_{0} -
480C_{0xx}\varphi_{x}\nonumber \\ & + 2400C_{0x}\varphi_{xx} +
1600D_{0y}\varphi_{x} - 300\varphi_{xx}A_{0}^3 +
1200\varphi_{xx}A_{0}C_{0} \nonumber \\ & - 2400\varphi_{xx}D_{1} -
39\varphi_{x}A_{0}^4 + 208\varphi_{x}A_{0}^2C_{0} -
400\varphi_{x}A_{0}D_{1} \nonumber \\ & + 400\varphi_{x}A_{1}D_{0} -
144\varphi_{x}C_{0}^2\Big].
 \label{pf1:15}
\end{eqnarray}
Since $\varphi=\varphi(x)$, we have $\beta_{y}=0$ yields
(\ref{eq:54})
$$
 \begin{array}{ll}
D_{1xy} =& \frac{3}{384}\Big[[(3A_{0}A_{1} - 4C_{1})A_{0}^2 +
16(2A_{1}D_{1} + C_{0}C_{1}) - 16(A_{1}C_{0} - D_{2})A_{0}]A_{0}
\\ & - 32[4(C_{1}D_{1} - 2C_{2}D_{0} + C_{0}D_{2}) + (3A_{1}D_{0} -
C_{0}^2)A_{1}]\\ & - 96D_{1y}A_{0} + 384D_{0y}A_{1} + 1536D_{0yy} -
16(3A_{0}A_{1} - 4C_{1})C_{0x} \\ & + 12[(3A_{0}A_{1} - 4C_{1})A_{0}
- 4(A_{1}C_{0} - 4D_{2})]A_{0x}\Big].
\end{array}
$$
From  (\ref{pf1:07}) one can rewrite the representation for $C_{0}$
upon denoting  $\chi = \frac{\varphi_{xx}}{\varphi_x}$ leads to
(\ref{eq:55}) and the representations for $\psi_{yy}$ and
$\psi_{xy}$ in the equations (\ref{pf1:01}) and (\ref{pf1:02})
become   (\ref{eq:57}) and (\ref{eq:58}).  Rewriting the
representation for $\alpha$ from  (\ref{pf1:12}) in the form
\[
    \alpha = \frac{\Omega}{8\varphi_{x}^3},
\]
where
\[
 \label{pf1:21}
 \Omega = A_{0}^3 - 4A_{0}C_{0} + 8D_{1} - 8C_{0x}
+ 6A_{0x}A_{0} + 4A_{0xx},
\]
and thus $\beta$ of  (\ref{pf1:15}) becomes
\begin{eqnarray}
\nonumber
 & \beta =  (1600\varphi_{x}^4)^{-1}( - 144A_{0x}^2
- 72A_{0x}A_{0}^2 + 352A_{0x}C_{0} + 160C_{0xx} + 80C_{0x}A_{0}
\nonumber \\ & + 1600D_{0y} - 640D_{1x} + 80\Omega_{x} - 9A_{0}^4 +
88A_{0}^2C_{0} - 160A_{0}D_{1} - 30A_{0}\Omega \nonumber \\ &+
400A_{1}D_{0} - 300\chi\Omega - 144C_{0}^2). \nonumber
\end{eqnarray}
Finally, one obtains (\ref{pf1:14}) in the form
\begin{eqnarray}
\nonumber & 1600\psi_{xxxx} =  9600\psi_{xxx}\chi + 160\psi_{xx}( -
12A_{0x} - 3A_{0}^2 - 90\chi^2 + 8C_{0}) \nonumber \\ & +
40\psi_{x}(12A_{0x}A_{0} + 72A_{0x}\chi - 16C_{0x} + 3A_{0}^3 +
18A_{0}^2\chi  - 12A_{0}C_{0} \nonumber \\ & + 120\chi^3 - 48\chi
C_{0} + 24D_{1} - 8\Omega)  + \psi(144A_{0x}^2 + 72A_{0x}A_{0}^2 -
352A_{0x}C_{0} \nonumber \\ &- 160C_{0xx} - 80C_{0x}A_{0}
 - 1600D_{0y} + 640D_{1x} - 80\Omega_{x} + 9A_{0}^4
- 88A_{0}^2C_{0} \nonumber \\ & + 160A_{0}D_{1}  + 30A_{0}\Omega -
400A_{1}D_{0} + 300\chi\Omega + 144C_{0}^2) + 1600\psi_{y}D_{0}.
\nonumber
\end{eqnarray}
Hence we complete the proof of Theorem 1.

\subsection{Proof of Theorem 2}
 \label{profcand1}

In the case of   (\ref{eq:21}), the problem is formulated as
follows. Given the coefficients $F_i(x,y), H_i(x,y), J_i(x,y),
K_i(x,y)$ of  (\ref{eq:21}), find the necessary and sufficient
conditions for integrability of the over-determined system of
equations (\ref{eq:23})-(\ref{eq:42}) for the unknown functions
$\varphi (x, y)$ and $\psi (x,y).$ Recall that, according to our
notation, the following equations hold
\begin{equation}
 \label{pf2:01}
 \varphi_{x} = r\varphi_{y}, \quad
 \psi_{x}=\frac{\psi_{y}\varphi _{x}-\Delta}{\varphi_{y}},
 \end{equation}
 and
\[
 \label{pf2:02}
 \alpha_{x}=\frac{\varphi_{x}}{\varphi_{y}}\,\alpha_{y},
  \qquad \beta_{x}=\frac{\varphi_{x}}{\varphi_{y}}\,\beta_{y}.
\]
Let us simplify the expression (\ref{eq:23}) as follow
\begin{equation}
 \label{pf2:03}
\varphi_{yy}  =  \Big[(4 \Delta_{y} - F_{2} \Delta)
\varphi_{y}\Big]/(10 \Delta).
\end{equation}
Comparing the mixed derivative
$(\varphi_{x})_{yy}=(\varphi_{yy})_{x}$  one obtains
\begin{eqnarray}\label{pf2:04}
 \Delta_{xy} = & \Big[F_{2x} \Delta^2 - F_{2y} r \Delta^2 + 10
r_{yy} \Delta^2 + 4 r_{y} \Delta_{y} \Delta - r_{y} F_{2} \Delta^2
\nonumber \\&+ 4 \Delta_{x} \Delta_{y} + 4 \Delta_{yy} r \Delta - 4
\Delta_{y}^2 r\Big]/(4 \Delta).
\end{eqnarray}
Rewriting  (\ref{eq:24}) in the form
\[
\Delta_{x}=  (20 r_{y} \Delta + 4 \Delta_{y} r + F_{1} \Delta - 2
F_{2} r \Delta)/4.
\]
Forming the mixed derivative $\Delta_{xy}=(\Delta_x)_{y}$ one
arrives at (\ref{eq:64})
\[
r_{yy} =   - (F_{1y} - F_{2x} - F_{2y} r - r_{y} F_{2})/10.
\]
Then (\ref{eq:25})-(\ref{eq:29}) are written in the form of
 (\ref{eq:65})-(\ref{eq:69}), respectively. Furthermore, (\ref{eq:30}) becomes
\[
\Delta_{yy}  =    - (20 F_{2y} \Delta^2 - 48 \Delta_{y}^2 + 4
\Delta_{y} F_{2} \Delta + 7 F_{2}^2 \Delta^2 - 20 J_{4}
\Delta^2)/(40 \Delta).
\]
Now, consider the equation  $(\Delta_{yy})_{x}=(\Delta_x)_{yy}$ ,
one gets  (\ref{eq:70})
\[
 \begin{array}{ll}
 F_{1yy} =  &  - (F_{1y} F_{2} - 40 F_{2xy} - 16 F_{2x}
F_{2} + 20 F_{2yy} r + 40 F_{2y} r_{y} \\ & + 14 F_{2y} F_{2} r  +
20 J_{4x} - 20 J_{4y} r + 14 r_{y} F_{2}^2 - 40 r_{y} J_{4})/10.
\end{array}
\]
Thus equations (\ref{eq:31})-(\ref{eq:34}) yield
(\ref{eq:71})-(\ref{eq:74}), and from (\ref{eq:35}) one finds
\begin{eqnarray}
\psi_{yyyy} = &\Big[300 \psi_{yyy} \varphi_{y} \Delta^2 (4
\Delta_{y} - F_{2} \Delta) +5 \psi_{yy}  \varphi_{y} \Delta ( - 120
F_{2y} \Delta^2 - 144 \Delta_{y}^2 \nonumber \\& + 72 \Delta_{y}
F_{2} \Delta - 39 F_{2}^2 \Delta^2 + 80  J_{4} \Delta^2) +\psi_{y}
\varphi_{y} ( - 500 \varphi_{y}^3 \alpha \Delta^3 \nonumber \\& -150
F_{2yy} \Delta^3 + 360 F_{2y} \Delta_{y} \Delta^2 -165 F_{2y} F_{2}
\Delta^3 + 100 J_{4y} \Delta^3 \nonumber \\& + 96 \Delta_{y}^3 - 72
\Delta_{y}^2 F_{2} \Delta + 108 \Delta_{y} F_{2}^2 \Delta^2 -240
\Delta_{y} J_{4} \Delta^2 - 24 F_{2}^3 \Delta^3 \nonumber \\& + 60
F_{2} J_{4} \Delta^3) -500 \psi \varphi_{y}^5 \beta \Delta^3 +500
K_{7} \Delta^4\Big]/(500 \varphi_{y} \Delta^3).
 \label{pf2:17}
  \end{eqnarray}
One can determine $\alpha$ from  (\ref{eq:36}), as the following
\begin{equation}
 \label{pf2:18}
\alpha =  (4 F_{2yy} + 6 F_{2y} F_{2} - 8 J_{4y} + F_{2}^3 - 4 F_{2}
J_{4} - 8 K_{6} + 56 K_{7} r)/8 \varphi_{y}^3.
\end{equation}
Now the equation  $\alpha_{x} - r \alpha_{y} = 0$ leads to
(\ref{eq:75}). Furthermore, one considers  $(\psi_x)_{yyyy} =
(\psi_{yyyy})_{x}$ , yields
\begin{equation}
 \label{pf2:20}
 \begin{array}{ll}
\beta = & 120 \Delta_{y} ( - 4 F_{2yy} - 6 F_{2y} F_{2} + 8 J_{4y} -
F_{2}^3 + 4 F_{2} J_{4} + 8 K_{6} - 56 K_{7} r) \\& + \Delta (320
F_{2yyy} + 480 F_{2yy} F_{2} + 336 F_{2y}^2 + 168 F_{2y} F_{2}^2 +
32 F_{2y} J_{4} \\&- 480 J_{4yy}  - 240 J_{4y} F_{2} - 1600 K_{7x} +
1600 K_{7y} r - 400 F_{1} K_{7} \\&- 9 F_{2}^4 + 88 F_{2}^2 J_{4}  +
160 F_{2} K_{6} - 320 F_{2} K_{7} r - 144 J_{4}^2)/1600 \Delta
\varphi_{y}^4.
\end{array}
\end{equation}
The equation  $\beta_{x} - r \beta_{y} = 0$ leads to (\ref{eq:76}).
Therefore,  (\ref{eq:37})-(\ref{eq:42}) become
 (\ref{eq:77})-(\ref{eq:82}), respectively.

Let us turn now to the integrability problem. One can find all
fourth-order derivatives of the functions $\varphi$ and $\psi$ by
using  (\ref{pf2:01}), (\ref{pf2:03}) and (\ref{pf2:17}). So that
one obtains at (\ref{eq:83})-(\ref{eq:86}). Finally, the
coefficients $\alpha$ and $\beta$ of the resulting linear equations
(\ref{pf2:18}) and (\ref{pf2:20}) are given by
\[
\alpha =  \frac{\Theta}{8 \varphi_{y}^3},
\]
\[
 \begin{array}{ll}
\beta =  & ( - 144 F_{2y}^2 \Delta - 72 F_{2y} F_{2}^2 \Delta + 352
F_{2y} J_{4} \Delta + 160 J_{4yy} \Delta + 80 J_{4y} F_{2} \Delta
\\& + 640 K_{6y} \Delta - 1600 K_{7x} \Delta - 2880 K_{7y} r \Delta -
4480 r_{y} K_{7} \Delta + 80 \Theta_{y} \Delta \\&- 120 \Delta_{y}
\Theta - 400 F_{1} K_{7} \Delta - 9 F_{2}^4 \Delta + 88 F_{2}^2
J_{4} \Delta + 160 F_{2} K_{6} \Delta \\&- 320 F_{2} K_{7} r \Delta
- 144 J_{4}^2 \Delta)/(1600 \varphi_{y}^4 \Delta),
\end{array}
\]
where
\[
\Theta = (F_{2}^2 - 4 J_{4}) F_{2} - 8 (K_{6} - 7 K_{7} r) - 8
J_{4y} + 6 F_{2y} F_{2} + 4 F_{2yy}.
\]
Hence we complete the proof of Theorem 2.

\section{Illustration of the linearization theorems}
\label{ex}

\subsection{An example on Theorem 1}
{\bf Example 1}. Consider the nonlinear ordinary differential
equation
 \begin{equation}
 \label{ex:01}
 x^2 y(2y^{(4)} + y) + 8x^2 y'y''' + 16xyy'''+ 6x^2y''^2 +
48xy'y''+ 24yy''+
 24y'^2=0.
 \end{equation}
It is an equation of the form (\ref{eq:09}) with the coefficients
\begin{eqnarray}
   & A_1 =  \frac{4}{y}\,, ~  A_0 =  \frac{8}{x}\,, ~  B_0 =
\frac{3}{y}\,, ~
 C_2 = 0\,, ~ C_1 = \frac{24}{xy}\,, ~ C_0 =
 \frac{12}{x^2}\,,~\nonumber\\[1.ex]
 & D_4 = 0\,, ~ D_3 = 0\,,
 ~ D_2 = \frac{12}{x^2 y}\,,~  D_1 = 0\,, ~ D_0 =
 \frac{y}{2}\,\cdot
 \label{ex:02}
\end{eqnarray}
One can check that the coefficients (\ref{ex:02}) obey the
conditions (\ref{eq:45})-(\ref{eq:54}). Thus, the equation
(\ref{ex:01}) is linearizable. We have
\begin{equation}\label{ex:03}
8C_0  - 3A_0^2  - 12A_{0x}  = 0
\end{equation}
and  the equation (\ref{eq:55}) is written as
\[
    2\frac{{d\chi }} {{dx}} - \chi ^2  = 0.
\]
Let us take its simplest solution  $\chi=0$ .  Then invoking
(\ref{eq:56}), we let
\[
   \varphi  = x.
\]
Now the equations (\ref{eq:57})-(\ref{eq:58}) are written
\[
    \frac{{\psi _{yy} }} {{\psi _y }} = \frac{1} {y},\quad \frac{{\psi
_{xy} }} {{\psi _y }} = \frac{2} {x}
\]
and yield
\[
    \psi _y  = Kx^2 y ,\quad \quad K = const.
\]
Hence
\[
   \psi  = K\frac{{x^2 y^2 }} {2} + f\left( x \right).
\]
Since one can use any particular solution, we set  $ K =
2,\;\;f\left( x \right) = 0$  and take
\[
    \psi  = x^2 y^2 .
\]
Invoking (\ref{ex:03}) and noting that (\ref{eq:60}) yields $\Omega
= 0 $, one can readily verify that the function  $ \psi  = x^2 y^2 $
solves equation (\ref{eq:59}) as well. Hence, one obtains the
following transformations
\begin{equation}\label{ex:04}
    t = x,\quad u = x^2 y^2.
\end{equation}
Since  $ \Omega  = 0 $,  equations (\ref{eq:61}) and (\ref{eq:62})
give
\[
   \alpha  = 0,\quad \beta  = \frac{1} {{\varphi _x^4 }} = 1
\]
Hence, the equation (\ref{ex:01}) is mapped by the transformations
(\ref{ex:04}) to the linear equation
\[
u^{\left( 4 \right)}  + u = 0.
\]
{\bf Example 2}. The third-order member of the Riccati Hierarchy is
given by Euler et al. \cite{Eul} as
\begin{equation}\label{ex:05}
    y'''+4yy''+3y'^{2}+6y^{2}y'+4y^{4}=0.
\end{equation}
Applying \cite{Mel}, and \cite{Sund} one checks that equation cannot
be linearized by a point transformation or contact transformation or
generalized Sundman transformation. Under the Riccati transformation
$y = \frac{a \omega'}{\omega}$ the equation (\ref{ex:05}) becomes
\cite{Leach}
\begin{eqnarray}
 &\omega ^3 \omega ^{\left( 4 \right)}  + 4\left( {a  - 1}
\right)\omega ^2 \omega '\omega ''' + 3\left( {a  - 1} \right)\omega
^2 \omega ''^2 \nonumber \\&+ 6\left( {a  - 1} \right)\left( {a  -
2} \right)\omega \omega '^2 \omega ''
  + \left( {a  - 1} \right)\left( {a  - 2} \right)\left( {a  - 3} \right)\omega '^4  = 0.
  \label{ex:06}
\end{eqnarray}
It is an equation of the form (\ref{eq:09}) with the coefficients
\begin{equation}
\label{ex:07}
\begin{array}{cc}
&A_1 =  \frac{4(a-1)}{\omega}\,, ~  A_0 = 0 \,, ~  B_0 =
\frac{3(a-1)}{\omega}\,, ~
 \\&
 C_2 = \frac{6 (a^2 - 3 a + 2)}{\omega^2 } \,, ~C_1 = 0\,, ~
 C_0 = 0\,,~
 \\&
 D_4 =\frac{a^3 - 6 a^2 + 11 a - 6}{\omega^3}\,, ~ D_3 = 0\,,~ D_2 = 0\,,~  D_1 = 0\,, ~ D_0
 =0.
\end{array}
\end{equation}
One can verify that the coefficients (\ref{ex:07}) obey the
linearization conditions (\ref{eq:45})-(\ref{eq:54}). Furthermore,
\begin{equation}\label{ex:08}
    8C_0  - 3A_0^2  - 12A_{0x}  = 0
\end{equation}
and  the equation (\ref{eq:55}) is written as
\[
    2\frac{{d\chi }} {{dx}} - \chi ^2  = 0.
\]
We take its simplest solution  $\chi=0$ and obtain from
(\ref{eq:56}) the equation $\varphi'' = 0$, whence
\[
   \varphi  = x.
\]
Equations (\ref{eq:57}) and (\ref{eq:58}) have the form
\[
    \frac{{\psi _{\omega \omega } }} {{\psi _\omega  }} = \frac{{a - 1}}
{\omega },\quad \quad \psi _{x\omega }  = 0
\]
and yield
\[
    \psi _\omega  = K \omega^{(a-1)},\quad \quad K = const.
\]
 Hence
\[
   \psi  = K\frac{\omega^{a}}{a} + f\left( x \right).
\]
Since one can use any particular solution, we set  $ K =
a,\;\;f\left( x \right) = 0$  and take
\[
    \psi  = \omega^{a} .
\]
Invoking (\ref{ex:08}) and noting that (\ref{eq:60}) yields $\Omega
= 0 $ , one can readily verify that the function  $ \psi  =
\omega^{a} $ solves equation (\ref{eq:59}) as well. So that one
obtains the following transformations
\begin{equation}\label{ex:09}
    t = x,\quad u = \omega^{a}.
\end{equation}
Since  $ \Omega  = 0 $,  equations (\ref{eq:61}) and (\ref{eq:62})
gives
\[
   \alpha  = 0,\quad \beta  = 0.
\]
Hence, the equation (\ref{ex:06}) is mapped by the transformations
(\ref{ex:09}) to the linear equation
\[
u^{\left( 4 \right)}   = 0.
\]
{\bf Example 3}. Let us consider the Boussinesq equation
\begin{equation}\label{ex:10}
    u_{tt}  + uu_{xx}  + u_x ^2  + u_{xxxx}  = 0.
\end{equation}
Of particular interest among the solutions of the Boussinesq
equation are travelling wave solutions:
\[
    u(x,t)=H(x-Dt).
\]
Substituting the representation of a solution into (\ref{ex:10}),
one finds
\begin{equation}\label{ex:11}
    H^{(4)} + (H + D^2)H'' + H'^{2} = 0.
\end{equation}
It is an equation of the form (\ref{eq:09}) with the coefficients
\begin{equation}
\label{ex:12}
\begin{array}{cc}
  & A_1 =  0\,, ~  A_0 = 0 \,, ~  B_0 = 0\,, ~
 C_2 = 0 \,, ~C_1 = 0\,, ~
 C_0 = D^{2} + H\,,~
 \\&
 D_4 =0\,, ~ D_3 = 0\,,~ D_2 = 1\,,~  D_1 = 0\,, ~ D_0
 =0.
\end{array}
\end{equation}
Since the coefficients (\ref{ex:12}) do not satisfy the
linearization conditions (\ref{eq:49}), (\ref{eq:52}) and
(\ref{eq:54}), hence, the equation (\ref{ex:11}) is not
linearizable.\\

\noindent {\bf Example 4}. Consider the non-linear equation
\begin{equation}\label{ex:13}
    y^{\left( 4 \right)}  - \frac{{10}} {{y'}}y''y''' + \frac{1} {{y'^2
}}\left( {15y''^3  - xy'^7  - y'^6 } \right) = 0.
\end{equation}
It has the form (\ref{eq:21}) with the following coefficients:
\begin{eqnarray}
  & r = 0\,, ~ F_2 =  0\,, ~  F_1 = 0 \,, ~  F_0 = 0\,, ~
 H_2 = 0 \,, ~H_1 = 0\,, ~
 H_0 = 0\,,~
 \nonumber \\&
 J_4 =0\,, ~ J_3 = 0\,,~ J_2 = 0\,,~  J_1 = 0\,, ~ J_0
 =0\,,~  K_7 = -x\,,~ \nonumber \\& K_6 = -1\,,~   K_5 = 0\,, ~
 K_4 =0\,, ~ K_3 = 0\,,~ K_2 = 0\,,~  K_1 = 0\,, ~ K_0
 =0.
 \label{ex:14}
\end{eqnarray}
 Let us test the equation (\ref{ex:13}) for linearization by using Theorem
 2. It is manifest that the equations (\ref{eq:64})-(\ref{eq:82})
are satisfied by the coefficients (\ref{ex:14}). Thus, the equation
(\ref{ex:13}) is linearizable, and we can proceed further.

Let us take its simplest solution $\varphi=y$ and $\psi=x$ which
satisfy  the compatible system  of equations
(\ref{eq:83})-(\ref{eq:86}). So that one obtains the following
transformations
\begin{equation}\label{ex:15}
    t = y,\quad u = x.
\end{equation}
Since  $ \Theta  = 8 $  ,  equations (\ref{eq:88}) and (\ref{eq:89})
give
\[
   \alpha  = 1,\quad \beta  =  1.
\]
Hence, the equation (\ref{ex:13}) is mapped by the transformations
(\ref{ex:15}) to the linear equation
\[
u^{\left( 4 \right)}  + u' + u = 0.
\]


\end{document}